\documentclass[12pt]{amsart}
\usepackage{geometry} 
\usepackage{latexsym, amsmath,amssymb}
\usepackage{enumerate}
\usepackage{color}

\usepackage{hyperref}
\newtheorem{theorem}{Theorem}[subsection]
\newtheorem{corollary}[theorem]{Corollary}
\newtheorem{proposition}{Proposition}[subsection]

\newtheorem{definition}{Definition}[subsection]

\def\F{\mathbf{F}}

\def\P{\mathbf{P}}

\def\N{\mathbf{N}}

\def\A{\mathcal{A}}

\def\cD{\mathcal{D}}

\def\cL{\mathcal{L}}

\def\cQ{\mathcal{Q}}
\def\cR{\mathcal{R}}

\title[QUANTUM SUSY OPERADS]{QUANTUM SUSY OPERADS}
\author{N. C. Combe, Yu.~Manin,  M. Marcolli } 

\date{} 

\begin{document}

\begin{abstract}  In a recent paper, we described a lifting of
coordinate rings of groups, loops, quantum groups, etc. to the categoric setup of  {\it operads}.
In most examples of that paper, these rings are non--commutative.

Quantum physics of the XX--th century added one more, quite nontrivial degree of freedom:
coordinates might become {\it fermionic}. In their classical
version, the fermionic coordinates {\it anti}--commute, and the resulting rings are called 
{\it supersymmetric}, or SUSY, ones.

In this paper, we try to lift operads involving fermionic coordinates 
to  {\it quantum operads}. We have to restrict ourselves
by lifting operads of supersymmetric rings. We also show that $1D$ supersymmetric algebras
have an operad structure, and we analyze their symmetries, through their relation to Adinkra graphs, dessins and codes.
\end{abstract}

\setcounter{section}{-1}
\maketitle

\tableofcontents
 {\it keywords}: monoidal categories, operads, super--symmetry.
\section{Introduction}~\label{S:Intro}

The lifting of classical structures to quantum ones generally goes through an initial stage that we call 
{\it categorification}:
basic data and axioms of the structure we deal with become categories, functors,
and commutative diagrams.

\smallskip

For this reason---which may seem a purely technical one---this paper turns out to be on the cross--roads
of two vast and differently motivated research domains: the general survey of the first one
can be found in \cite{Ke01}; as for the second one, it can be found in \cite{KaSh06}. How disjoint they might
seem to researchers, a reader can guess by simply comparing the lists
of References, in these two surveys. The key references appear either in one
of them, or in another, but rarely in both.

\smallskip

We have chosen \cite{KaSh06} as our main source of basic information.

\smallskip

The following four papers focus on the central objects of our current  study:

\smallskip
\begin{list}{--}{}
\item Operads and generalised operads: \cite{BoMa08};

\item Operads of moduli spaces of algebraic curves: \cite{GeKa98};

\item Quantum operads of moduli spaces of algebraic curves: \cite{CoMaMar22};

\item Generalised operads of moduli superspaces of stable SUSY curves: \cite{KeMaWu22}.

\end{list}
\smallskip

Briefly, the goal of this paper is to present a combination of quantisation and supersymmetry.
For further preparatory reading, we can recommend to a reader the following references
\cite{MaPe89,FeKaPo20,HoKrStT11,WuYau06}.

\smallskip

The sequence, in which we explain and briefly discuss the potential setup
of quantum SUSY modular operads, is motivated by the presentation
of quantum operads, not involving fermionic coordinates, in
\cite{CoMaMar22,HoKrStT11}, and the survey~\cite{Sm16}. See also
the introductory parts of \cite{BoMa08} and \cite{MaVa20}.

At the end of this paper, we discuss another source of quantum SUSY operads: $\bold{F}_2$-codes and Adinkras and investigate their hidden symmetries.

\bigskip
 
 \section{Categorifications and enrichments: the simplest examples}

 All our definitions and constructions below will refer to {\it sets from
 a fixed universe} $\mathcal{U}$: see \cite{KaSh06}, Sec. 1.1, pp. 10--11. In particular, we will use basic notions of the language of categories
 as they are presented in \cite{KaSh06}, in Sec. 1.2, and further on.
 
 \medskip
 
Somehow, this foundational book does not stress the importance (also historical importance)
 of the language of {\it structures} \`a la Bourbaki. So, we start this article
 with an attempt to formalise and generalise interrelationships between these two
 notions. 
 
 \smallskip
 
 On the one hand, we discuss {\it  transitions from structures to
 categories}, which we will call {\it categorifications}.  On the other hand, we discuss {\it transitions from a category
 to ``structured objects of category"}, which then form another category, being an {\it enrichment}
 of the initial one (or of its part).
 
 \smallskip
 
Let us start with basic examples.
 
 \medskip
 
\subsection{Monoids and monoidal categories} 
A {\it monoid} $M$, as it is defined on p.13 of \cite{KaSh06}, is a set, endowed with binary multiplication
 $m:\, M\times M \to M$, satisfying the associativity  law. It might also have left, resp. right, 
 identities. If there is one element $\mathbf{1}$, which is simultaneously left
 and right, it is called simply identity.
 
 \smallskip
 
 Recall, that for two $\mathcal{U}$--sets $X,Y$, their direct product $X\times Y$ is the 
$\mathcal{U}$--set,
 whose elements are ordered pairs $(x,y)$ with $x\in X,\, y\in Y$.
 If we write $(xy)$ in place of $m(x,y)$, the associativity condition
 can be written in the standard way: $((xy)z)= (x(yz))$. (Exterior brackets
 are usually omitted, but they become necessary for a linear ordering 
 of iterated operation).
 
 \smallskip
 
 Categorification of monoids leads to the definition of a monoidal
 category (cf. \cite{Sm16},  Sec. 2.2, 2.3,  and \cite{KaSh06}, Sec. 4.2, 
 where monoidal categories are called tensor ones).
 
 \medskip
 
\begin{definition}  A monoidal category $\mathcal{T}$  is a category 
 whose set of objects  is endowed with associative multiplication
 and identity, functorial with respect to morphisms in this category.
 \end{definition}
 
 \smallskip
 
 More precisely, a monoidal structure on $\mathcal{T}$ is given by the following
 data, which we call {\it structural functors}: binary
 multiplication $\otimes$ of objects, with identity object $\mathbf{1}$ and natural isomorphisms
\[
\alpha_{A,B;C}: (A\otimes B)\otimes C  \to A\otimes (B\otimes C), \quad
\rho_A: A\otimes \mathbf{1} \to A, \quad
\lambda_A: \mathbf{1}\otimes A \to A .
\]
See the diagram (4.2.1) on p.96 of \cite{KaSh06} for a detailed
description of the categorification of associativity law.

\smallskip
We refer to this set of data as {\it structural functors}.

\medskip

\subsection{Commutative monoids and symmetric monoidal categories} 
A monoid $(M,m)$ is commutative, if $m(x,y)=m(y,x)$ for all $x,y \in M$.

\smallskip

Respectively, the list of structural functors in the definition of a symmetric monoidal category
 includes additional {\it twist} isomorphisms
$\tau_{A,B}: A\otimes B \to B\otimes A$, with $\tau_{A,B}\tau_{B,A} = \mathrm{id}_{A\otimes B}$,
satisfying axioms of compatibility with associativity, expressed by several commutative diagrams.

\medskip
 
\subsection{Magmas, comagmas, bimagmas in symmetric monoidal categories} 

The basic data for a magma: an object $A$ with multiplication morphism $\nabla : A\otimes A \to A$.

\smallskip

The basic data for a comagma: an object $A$ with comultiplication morphism $\Delta : A\to A\otimes A$.

\smallskip 

The basic data for a bimagma: a triple $(A,\nabla , \Delta)$ as above such that the ``bimagma diagram''
 ((2.4)~\cite{Sm16}, p.49) commutes.

\smallskip

(Co, bi)--magmas in a symmetric monoidal category $\mathbf{V}$ are themselves objects
of respective categories. Morphisms of them are those morphisms in $\mathbf{V}$, (\cite{Sm16}, Def. 2.3),
which are compatible with the respective basic data.

\medskip

\subsection{Commutative/cocommutative  magmas and comagmas in symmetric monoidal categories} 
(\cite{Sm16}, Def.2.3). Basically, these properties mean the compatibility with corresponding
structural functors.

\smallskip

\subsection{Monoids, comonoids, bimonoids, and Hopf algebras in symmetric mono\-idal categories}  (\cite{Sm16}, Def. 2.7).
They are essentially (co, bi)--magmas with additional (co,bi)--associativity restrictions.

\medskip

\subsection{Quantum quasigroups} (\cite{Sm16}, Sec. 3.1). 
{\it A quantum quasigroup} $(A, \nabla , \Delta)$
is  a bimagma, for which both {\it left composite} and {\it right composite} morphisms
are invertible:
\[
A\otimes A  \xrightarrow{\Delta\otimes \mathrm{id}_A} A\otimes A \otimes A \xrightarrow{\mathrm{id}_A\otimes \nabla} A\otimes A,
\]

\[
A\otimes A \xrightarrow{\mathrm{id}_A\otimes \Delta }  A\otimes A \otimes A  \xrightarrow{
\nabla \otimes \mathrm{id}_A} A\otimes A .
\]

\smallskip

These morphisms are sometimes called {\it fusion operators} or {\it Galois operators}.

\smallskip

\subsection{Quantum loops} 
{\it A quantum loop} in $\mathbf{V}$ is a biunital bigmagma $(A,\nabla ,\Delta, \eta, \varepsilon )$
such that $(A, \nabla ,\Delta )$ is a quantum quasigroup. (The notation $(\eta , \varepsilon )$
is introduced in [Sm16], Sec. 2.1).

\smallskip

\subsection{Functoriality} (\cite{Sm16}, Prop. 3.4). Any symmetric monoidal functor 
\[
F: (\mathbf{V}, \otimes , \mathbf{1}_{\mathbf{V}} ) \to ( \mathbf{W}, \otimes , \mathbf{1}_{\mathbf{W}})
\]
sends quantum quasigroups (resp. quantum loops) in $\mathbf{V}$ to quantum quasigroups
(resp. quantum loops) in $\mathbf{W}$.

\smallskip

\subsection{Magmas etc. in   the categories of sets with direct product}
 According
to \cite{Sm16}, beginning of Sec. 3.3, in such categories comultiplication in a counital
comagma is always the respective diagonal embedding.
As a corollary, we see that quantum loops and counital quantum quasigroups in
such caregories are cocommutative and coassociative.

\smallskip

As a result, we see, that in such a category, counital quantum quasigroups are equivalent
 to classical quasigroups, and quantum loops are equivalent to classical loops
 (\cite{Sm16}, Prop. 3.11).

\bigskip

\section{Monoidal categories of operads} 

\medskip

\subsection{Graphs and their categories} 
Our basic definition of graphs as quadruples $(F,V,\partial, j)$ and their categories is explained in \cite{BoMa08}, Sec. 1.1, p.251.  There $F$, resp. $V$, are called the sets of {\it flags}, resp. {\it vertices},
 and structure maps $\partial$, resp. $j$ are called {\it boundary maps}, resp. {\it involutions}.
 Usually one flag is a pair consisting of flag as such, and a {\it label},
 that should be defined separately. 
 
 {\it Geometric realization} of a graph
 is the quotient set of the disjoint union of semi--intervals $(0, 1/2]$ labeled with flags
 of this graph, modulo equivalence relation, in which $0$--points of a flag
 is glued to $1/2$ of another flag, if these flags are related by the boundary
 relation, or structure involution.
 
 \smallskip
 
 Depending on the context and/or type of labelling of $\tau$, elements of $F_{\tau}$
 might be called {\it flags, leaves, tails} ... In the study of magmatic operad
 (\cite{ChCorGi19}) and the relevant binary trees, vertices of the relevant corollas
 are called {\it nodes}, non--root flags are called {\it left child, right child} etc.
 We will try to attach all such ``heteronyms'' to our basic terminology of \cite{BoMa08}.
 
 \smallskip
Below the most typical labeling of our graphs will be (see details in \cite{BoMa08},
Sec. 1.3.2 a) and 1.3.2 e), pp. 257--259):

\vspace{3pt}

\begin{enumerate}[(i)]
\item {\it  Orientation.}\\
\item {\it Cyclic labeling.}
\end{enumerate}

\vspace{3pt}

To give an orientation and cyclic labeling of corolla is essentially the same as
to define it as a {\it planar} graph: corolla, embedded into an oriented real
affine plane, with labeling compatible with its orientation.

\smallskip

Graphs endowed with various labelings form categories, upon which
the operation of disjoint union $\sqcup$ defines a monoidal structure: see \cite{BoMa08},
Sec. 1.2.4, pp. 254--255. Our central objects of study are initially defined
only for connected graphs. Therefore, introducing this monoidal product,
we must first take care of ``empty'' (or partially empty) graphs and
explain details of their functoriality. The paper \cite{BoMa08} is interspersed
with subsections directly or indirectly motivated by this necessity.

\smallskip

For the purposes of this paper, the most important graphs are labelled {\it trees} and
{\it forests} -- disjoint unions of trees, forming {\it ``selva selvaggia e aspra e forte''}.

\medskip

\subsection{Operads and categories of operads}  (See \cite{BoMa08}, Sec. 1.6, p.262). 
We recall here the first definition of operads in \cite{BoMa08}, 1.6 (I), and morphisms of operads as in \cite{BoMa08}, Sec. 1.6.1.

\smallskip

First of all, we fix a symmetric monoidal category of labelled graphs $\Gamma$
with disjoint union as the monoidal structure, and a symmetric monoidal
{\it ground category} $(\mathcal{G}, \otimes)$, satisfying a part of conditions 1.4 a) -- f) in \cite{BoMa08}, p.259.

\vspace{3pt}

\begin{enumerate}[(i)]
\item  {\it An operad is  a tensor functor between two monoidal categories
$A: (\Gamma , \sqcup )\to (\mathcal{G},\otimes )$ that sends any grafting morphism to
an isomorphism.}
\\
\item  {\it A morphism between two operads is a functor morphism.}  
\end{enumerate}

\vspace{5pt}

Denote this category of operads by $\Gamma\mathcal{G}OPER$.

\medskip

\subsection{Operads and collections as  symmetric monoidal categories}
 Following \cite{BoMa08}, Sec. 1.8, we
will introduce now the monoidal ``white product'' of two operads 
$A,B : (\Gamma ,\sqcup ) \to (\mathcal{G}, \otimes)$ by
the formula
\[
A\circ B (\sigma ) := A(\sigma) \otimes B(\sigma)
\]
extended to morphisms in a straightforward way.

\smallskip

Clearly, $(\Gamma\mathcal{G}OPER, \circ)$ is a symmetric monoidal category.

\smallskip

An important related notion is that of {\it collection}. Starting with $\Gamma$ as above,
denote by $\Gamma COR$ its subcategory, whose objects are corollas in $\Gamma$,
and morphisms between them are isomorphisms.

\smallskip

Combining  it with the ground category $(\mathcal{G}, \otimes)$ as obove, we can introduce
the category $\Gamma\mathcal{G} COLL$ of $\Gamma G$--collections: its
objects are functors $A_1 :  \Gamma COR \to \mathcal{G}$, and morphisms are 
natural transformations between these functors.

\smallskip

The restriction of white product $\circ$ to $\Gamma\mathcal{G} COLL$ defines on it the structure
of symmetric monoidal category. If $(\mathcal{G}, \otimes )$ has an identity object $\mathbf{1}$,
then the collection $\mathbf{1}_{coll}$ sending each corolla to $\mathbf{1}$ and each
isomorphism of corollas to the identical isomorphism of $\mathbf{1}$, is the 
identity collection.

\medskip

The restriction of white product $\circ$ to $\Gamma\mathcal{G} COLL$ defines on it the structure
of symmetric monoidal category. If $(\mathcal{G}, \otimes )$ has an identity object $\mathbf{1}$,
then the collection $\mathbf{1}_{coll}$ sending each corolla to $\mathbf{1}$ and each
isomorphism of corollas to the identical isomorphism of $\mathbf{1}$, is the 
identity collection.

\medskip

\subsection{Operads as monoids}
 We briefly describe here a construction by
B. Vallette~\cite{Val04}, reproduced in~\cite{BoMa08}, Appendix A, Subsection 5.

\smallskip

We will have to use here a stronger labeling of graphs in $\Gamma$
than just orientation. Besides orientation, connected objects of $\Gamma$ 
must admit a continuous real--valued function such that it decreases
whenever one moves in the direction of orientation along each flag.
Such graphs are called {\it directed ones}
(see~\cite{BoMa08}, Sec. 1.3.2 b). 
\smallskip

A graph $\tau$ is called {\it two--level graph}, if it is oriented, and if there
exists a partition of its vertices $V_{\tau} = V^1_{\tau} \sqcup V^2_{\tau}$
with the following properties:

\vspace{3pt}

\begin{enumerate}[(i)]
\item  {\it Tails at $V^1_{\tau}$ are all inputs of $\tau$, and tails
at $V^2_{\tau}$ are all outputs of $\tau$.}
\\
\item  {\it All  edges in  $E_{\tau}$ go from $V^1_{\tau}$ to $V^2_{\tau}$.}
\end{enumerate}

\vspace{5pt}

For any two $\Gamma G$--collections $A^1, A^2$ define their product as
\[
(A^2 \boxtimes_c A^1) (\sigma ) := \mathrm{colim} 
(\otimes_{v\in V^1_{\tau}} A^1(\tau_v)) \otimes (\otimes_{v\in V^2_{\tau}} A^2(\tau_v)).
\]

Here colim is taken over the category of morphisms from two level graphs to $\sigma$.

\smallskip

\begin{theorem} The product $\boxtimes_c$ is a monoidal structure on collections,
and operads are monoids in the respective monoidal category.
\end{theorem}
\smallskip

See \cite{Val04} and \cite{BoMa08}.

\medskip

\noindent  2.4.2 \noindent {\bf Freely generated operads} For any $\Gamma\mathcal{G}$--collection $A_1$ one can define
another collection $\mathcal{F}(A_1)$ together with a canonical structure of operad on it, and for any operad
$A$ each morphism of collections $A_1 \to A$ extends to a morphism of operads
$f_A: \mathcal{F}(A_1) \to A$. 

\smallskip

We can imagine $\mathcal{F}(A_1)$ as the operad freely generated by the collection $A_1$.

\smallskip

\subsection{Comonoids in operadic setup}
 We will now introduce a category $OP$ of operads
given {\it together with their presentations} (\cite{BoMa08}, Sec. 2.4). We start with $\Gamma$ and $\mathcal{G}$
as above.

\smallskip

One object of $OP$ is a family $(A, A_1, i_A)$, where $A$ is a $\Gamma\mathcal{G}$--operad,
$A_1$ is a $\Gamma\mathcal{G}$--collection, such that  $f_A: \mathcal{F}(A_1) \to A$ is surjective.

\smallskip

Define on $OP$ a product $\odot$ by the formula
\[
(A, A_1, i_A) \odot (B, B_1, i_B) = (C, C_1, i_C) ,
\]
in which $C_1 := A_1\circ B_1$ (cf. 2.3 above), $C$:= the minimal suboperad,
containing the image $(i_A\circ i_B)(A_1\circ B_1) \subset A\circ B$,
and $i_C$ is the restriction of $I_A\circ i_B$ on $A_1\circ B_1$.

\smallskip

\begin{theorem} (See \cite{BoMa08}, Sec. 2.4).
\ 
\begin{enumerate}[(i)]
\item The product $\odot$ defines on $OP$ a structure
of symmetric monoidal category.\\

\item The category $OP$ is endowed with the functor of inner cohomomorphisms
\[
\underline{cohom}_{OP} : OP^{op} \times OP \to OP
\]
so that we can identify, functorially with respect to  all arguments,
\[
\mathrm{Hom}_{OP} (A, C\odot B) = \mathrm{Hom}_{OP} (\underline{cohom}_{OP} (A,B), C)
\]

\item Therefore, one can define canonical coassociative comultiplication
morphisms
\[
\Delta_{A,B,C} : 
\underline{cohom}_{OP} (A, C) \to \underline{cohom}_{OP}(A,B) \odot \underline{cohom}_{OP} (B,C) .
\]
\end{enumerate}
\end{theorem}
\smallskip

\begin{corollary}   For any $A$, the coendomorphism operad
\[
\underline{coend}_{OP} A := \underline{cohom}_{OP} (A,A)
\]
is a comagma in the sense of 1.3 above.
\end{corollary}

\smallskip

\subsection{The magmatic operad}  (See \cite{ChCorGi19}). 
Below we give a brief survey
of some definitions and results from \cite{ChCorGi19}, sometimes slightly changing terminology and
notation.

\smallskip

Here objects of our basic symmetric monoidal category $(\Gamma , \sqcup )$ will be
disjoint unions of oriented trees with the following additional labeling: {\it for each tree, its outcoming flags (or leaves)
are cyclically ordered}. Corollas in it are one--vertex graphs with one root and at least
two leaves. Connected objects can be obtained from a union of disjoint corollas
by grafting each root of a corolla to one of leaves of another corolla. Morphisms
are compatible with labeling.

\smallskip

An algebra over magmatic operad is a family $(\A,*)$ consisting of a set $\A$ with binary composition law
$* : \A\times\A \to \A$.

\smallskip

Thus, corollas in the magmatic category correspond to products 
\[
(x_1*((x_2)* \dots (...(x_n)))...),
\]
and generally, connected graphs in it correspond to monomials of generic arguments with all possible
arrangements of brackets.

\smallskip

\subsection{Quasigroup monomials and planar trees} 
Monomials that can be obtained
by iteration of binary multiplication $*$ as in (0.1) correspond to  {\it planar} trees:
see 2.1 above for discussion of planar corollas.
Below, discussing  quasigroups in general, and Moufang loops in particular,
we will consider connected planar trees and quasigroup monomials as encoding each other
in this way.

\bigskip

\section{Categories of quadratic data and their relationships to supersymmetry}

The exposition in this section starts with structures, considered  in \cite{Ma88,MaVa20},
and proceeds to their categorifications and (partial) enrichments.

\smallskip

\subsection{Quadratic data}
 According to Sec. 2 of \cite{MaVa20}, one object
of this category is a pair $(V,R)$, consisting of a finite--dimensional  $\mathbf{Z}$--graded
vector space $V$ (over a field of characteristic $\neq 2$), and a subspace $R\subset V^{\otimes 2}$.

\smallskip

A morphism $f:  (V,R) \to (W,S)$ is a morphism of graded vector spaces $f: V\to W$,
such that $f^{\otimes 2} (R) \subset S$.

\medskip

\subsubsection{Quadratic data in a symmetric monoidal category} 
 Let $\mathcal{T}$ be a symmetric monoidal category.

\smallskip

The respective enrichment of the category of quadratic data has---as its objects---
pairs $(V,R)$, where now $V$ is an object of $\mathcal{T}$, and
$R$ is a subobject of $V^{\otimes 2}$. A morphism $f: (V,R) \to (W,S)$,
is a morphism $f: V \to W$ in $\mathcal{T}$,
satisfying the same restriction, as above, in 3.1.

\medskip

\subsection{A SUSY categorification of symmetric monoidal categories} The {\it data} is defined  
as above, and equipped with a functor {\it ``sign change"} $c: \mathbf{V} \to \mathbf{V}$
such that $c\circ c = Id_{\mathbf{V}}$. 

\smallskip

Moreover, we require $c(\mathbf{1}) = \mathbf{1}$.

\smallskip

For example, in the category of $\mathbf{Z}_2$--graded vector spaces,
$c$ acts as sign change on the oddly graded vectors.

\smallskip

This allows us to imitate constructions, using $\mathbf{Z}_2$--gradings,
when we pass from a commutative setup to the supercommutative one.

\medskip

\subsection{SUSY geometry}
 The main constructions, upon which we
will focus further,  are based upon systems of notions
developed in \cite{De87} and then specialised in the domain
of stable SUSY families of curves. We omit this
general setup. An interested reader can turn to other references:
see \cite{BrHePo20}.

\bigskip

\section{Operads of moduli spaces of SUSY curves}

\medskip

This section starts with a brief exposition of the recent work \cite{KeMaWu22}, whose main goal
is a construction of the extension of the classical modular operad, taking into account
enrichments of curves to supercurves. 

\smallskip

More precisely, stable degenerations in families of classical
curves, describing strata in moduli spaces, involve gluing pairs of smooth points into double points,
or {\it nodes}.

\smallskip

When we pass from curves to supercurves, the stability restrictions
become considerably stricter, and as was shown in \cite{De87}, only
two types of nodes do not break superstability: {\it Neveu--Schwarz} ones and {\it Ramond} ones.

\smallskip

Finally,  as we proceed to quantisation, additional difficulties arise:
modular operads in the style of \cite{GeKa98} can be constructed {\it only}
from superspaces parametrising SUSY {\it curves  of genus zero}.
The main reason for it can be intuitively described as follows:
whenever we work with the usual operadic maps for stacks of supermoduli of
higher genus, the usual operadic morphisms lose information about the relevant spinor bundles,
and therefore cannot combine to a full operad. For details,
see \cite{KeMaWu22}, Sec. 4.1.

\smallskip

The existence of the quantum genus zero SUSY operad, restricted to
Neveu--Schwarz and Ramond punctures, is the principal new result of this paper.

\medskip

\subsection{Category of graphs, encoding SUSY curves} 
Start with choosing
of $k_{NS}\in \bold{Z}_{\ge 0}$ and $k_R \in 2\bold{Z}_{\ge 0}$.

\smallskip

According to the
Sec. 2.3 of \cite{KeMaWu22},  {\it one object} of this
category is {\it a $(k_{NS}, k_R)$--SUSY (or simply a SUSY--graph) $\tau = (F_{\tau}, V_{\tau}, \partial_{\tau}, j_{\tau})$,
endowed with}

\vspace{3pt}
\begin{enumerate}[(i)]
\item {\it a genus labeling} $g_{\tau} :  V_{\tau} \to \bold{Z}_{\ge 0}$;\\
\item {\it a puncture colouring}  $c_{\tau} : F_{\tau} \to \{NS , R \}$,
such that each vertex $v \in V_{\tau}$ the number of adjacent to this vertex flags with color
$R$ is even;\\
\item {\it Two separate labelings} of $NS$--tails and $R$--tails, that
are bijections
\[
l_{\tau ,NS} :  \{1, \dots , k_{NS} \} \to T_{\tau, NS} := F_{\tau,NS}\cap T_{\tau} ,
\]
\[
\ l_{\tau, R} : \  \{1,  \dots , k_R\}\  \to T_{\tau, R}\  := F_{\tau, R} \cap T_{\tau}.
\]
\end{enumerate}
When these structures are chosen, then the sets of flags, tails, and edges,
adjacent to each vertex $v\in V_{\tau}$, acquire colours $NS$ ({\it Neveu--Schwarz})
or $R$ ({\it Ramond}), and are called respectively.

\smallskip

{\it A morphism} in this category is a morphism of respective graphs in the sense
of \cite{BoMa08}, compatible with the genus labelings and colourings of flags.
Some tails of the same colour can be grafted, and some pairs of the tails of the same
colour can be virtually contracted.

\medskip

\subsection{SUSY curves with punctures} 
 The notions of SUSY curves, their families,
morphisms, etc. are natural extensions of the respective notions in the formalism
of schemes, starting with a replacement of structure sheaves of commutative rings of schemes
by structure sheaves of $\mathbf{Z}_2$ graded supercommutative rings of superschemes.

\smallskip

Skipping these foundational preparations, we pass to the description of
relevant notions of SUSY curves with punctures, their families, and
notions of stability, from \cite{De87} and \cite{FeKaPo20}:

\smallskip

\begin {definition}  \label{D:4.2.1}
A family of SUSY curves with $k_{NS}$ Neveu--Schwarz
punctures and $k_R$ Ramond punctures over the base $B$ consists of the following data:

\vspace{3pt}
\begin{enumerate}[(i)]
\item A smooth proper morphism of superschemes $\pi : M \to B$, whose generic fibres
have relative dimension $1|1$.

\smallskip

\item  A sequence of sections $s_i :  B \to M$, $i = 1, \dots k_{NS}$, such that
on each fibre of $\pi$ the reductions of $s_i$ and $s_j$ are different for $i\ne j$.

These reductions are called Neveu--Schwarz punctures.

\smallskip

\item   A sequence of components $r_j $, $j = 1, \dots , k_R$
of an unramified effective Cartier divisor $\mathcal{R}$ of codimension $0|1$ 
of  degree $k_R$

These components are called Ramond punctures.

\smallskip

\item  The line bundle $\mathcal{D}$, a subbundle of the tangent bundle $\mathcal{T}_M$
of rank $0|1$, such that the commutator of vector fields induces an isomorphism
\[
\mathcal{D} \otimes \mathcal{D} \to (\mathcal{T}_M  / \mathcal{D} ) (- \mathcal{R}) .
\]
\end{enumerate}
\end{definition}
\medskip

\subsection{Stability} \label{S:4.3}
Consider a family of SUSY curves with punctures 
over $B$ as above. It is called {\it (super)stable}, if it
satisfies the following restrictions:

\vspace{3pt}
\begin{enumerate}[(i)]
\item {\it $M$ is a proper, flat and relatively Cohen--Macaulay superscheme over $B$}.\\\item {\it $M$ contains an open fibrewise dense subset $U$ containing
all sections $s_i$ and $r_j$ and such, that $U/B$ is smooth
of relative dimension $1|1$.}\\
\item {\it The reduction $M_{red} \to B_{red}$ is a stable family
of marked curves.}
\end{enumerate}
\bigskip

\subsection{Dual graphs of stable SUSY families} 
The dual graph $\tau$,
associated with a stable family of SUSY curves (as in Def.~\ref{D:4.2.1} and  Sec.~\ref{S:4.3}
such that $B_{red}$ is a point) is defined as follows:

\vspace{3pt}
\begin{enumerate}[(i)]

\item {\it The set of its vertices $V_{\tau}$ is identified with the set of irreducible
components of M}.

\smallskip

\item {\it Flags of $\tau$ with boundary $v$ are identified with special points
on the respective irreducible component of $M$}.

\smallskip

\item  {\it  The involution $j_{\tau}$ changes places of halves of edges
with one boundary $v$.} In particular, it distinguishes punctures: those
belonging to different components of a stable curve and those
corresponding to double points of one component.\\
\item {\it The genus labeling marks irreducible components by their
genera.}\\
\item {\it The type of a puncture breaks the whole set of flags into two
disjoint subsets: $F_{\tau} =F_{\tau , NS} \cup F_{\tau , R}$.}\\
\item {\it The labelings $l_{\tau , NS}$, resp. $l_{\tau , R}$, mark flags,
corresponding to respective punctures}.
\end{enumerate}
\bigskip

\subsection{Moduli stacks of stable SUSY families}
 Generally, the functor
of stable SUSY families of curves of genus $g$  with a fixed dual graph is represented 
by a smooth and proper Deligne--Mumford superstack
$\overline{\mathcal{M}}_{g, k_{NS} , k_R}$ : this was proved in 
\cite{FeKaPo20}. 

\smallskip

The families of curves with at least one Neveu--Schwarz node,
resp. at least one Ramond node,
are represented by the boundary Cartier divisor $\Delta_{NS}$,
resp. $\Delta_R$.

\smallskip

A fundamental role in the construction of operadic compositions for 
moduli stacks is played by glueing punctures of stable SUSY curves.

\smallskip

An attempt to lift all operadic morphisms from the classical to the SUSY setup
was made in \cite{KeMaWu22}, Sec. 4.

\smallskip

In this final part of our study we draw attention of the reader to the fact,
that if we want to have these liftings and canonical morphisms
among them to be unique, so that we get an actual operad
on the level of SUSY moduli spaces, we have to restrict ourselves
by genus zero and Neveu--Schwarz and Ramond punctures.

\medskip

\subsection{SUSY modular operad breaking down}
 The breaking of higher genus morphisms is explained in \cite{KeMaWu22}, Sec. 4.1: any kind of glueing
requires looking at higher general components: see equation (4.1.2), (4.1.3), (4.1.4).

\smallskip

However, if we consider only genus zero components and glueing of
Neveu--Schwarz, resp. Ramond, punctures, the same diagrams
become united into one operad.

\medskip 

\section{Operadic structure of supersymmetry algebras}

\bigskip

In Section~5 of \cite{CoMaMar22}, it was shown that certain classes of classical and quantum
codes carry the structure of an algebra over a version of the little squares operad.
In this section we show that the set of representations of $1D$ supersymmetry algebras 
carries the structure of an operad.

\smallskip

We first recall some general facts about $1D$ supersymmetry algebras and their
description in terms of Adinkra graphs and linear codes. We then use the description
in terms of codes to introduce the operadic structure. 

\smallskip
\subsection{Supersymmetry algebras and codes}\label{S:5.1}
In the setting of supersymmetric quantum mechanics, with a $1$-dimensional space-time
with time coordinate $t$ and a zero-dimensional space, the $(1|N)$ Superalgebras are
generated by operators  $Q_1, \, Q_2,\, \cdots,\, Q_N$, which give the supersymmetry generators, 
and $H=i\partial_t$, subject to the commutation and anticommutation relations: 
\begin{equation}\label{E:5.1}
 [Q_k, H] = 0 \ \ \ \text{ and } \ \ \ 	\{Q_k, Q_\ell\} = 2 \delta_{k\ell} H,  
 \end{equation}
 where $\delta_{k\ell}$  stands for the Kronecker delta.
 
Representations of these algebras on bosonic and fermionic fields, $\{\phi_1, \ldots, \phi_m \}$
and $\{\psi_1, \dots, \psi_m \}$, respectively, are of the following form:
\begin{equation}\label{E:5.2} 
Q_k \phi_a = c\, \partial_{t}^\lambda\,  \psi_b, 
\end{equation}
\begin{equation}\label{E:5.3}
	Q_k \psi_b = \frac{i}{c}\, \partial_{t}^{1 -\lambda}\, \phi_a \, ,
	\end{equation}
with $c \in \{-1,1\}$ and $\lambda \in \{0,1\}$. 

\smallskip

A graphical method for classifying these  supersymmetry algebras representations was introduced by Faux and Gates, (see \cite{FaGa05}), in terms of decorated bipartite graphs with additional structure, called {\it Adinkras}.
The geometry of Adinkras, their relation to a class of Grothendieck's dessins d'enfant,
and their classification in terms of linear codes were developed in  \cite{DFGHIL07}, and \cite{DIKLM15,DIKLM16}.

\smallskip

Summarizing the correspondence between $1D$ SUSY algebras and Adinkras,
we have the following setting, \cite{FaGa05}.

\smallskip

Let $A$ be a finite graph with no looping edges and no parallel edges. Let $V(A)$ and $E(A)$
be the sets of vertices and edges. Such a graph is an {\it  $N$-dimensional chromotopology} if
it satisfies the following properties:

\vspace{3pt}
\begin{enumerate}[(a)]
\item   $A$ is $N$-regular and bipartite, with vertices in the bipartition colored white and black respectively;

\item  The edges in $E(A)$ are colored by $N$ colors, labelled by $\{1, 2, \cdots, N\}$,
with every vertex incident to exactly one edge of each color;
\smallskip

\item  For any pair $i\neq j$ of colors the edges in $E(A)$ labelled with 
colors $i$ and $j$ form a disjoint union of $4$-cycles.
\end{enumerate}

\smallskip

A {\it ranking} of the graph $A$ is a partial ordering of $V(A)$ determined by
a function $h: V(A) \rightarrow \mathbf Z$. A {\it dashing} of the graph $A$
is a function $d: E(A) \rightarrow \mathbf F_2$ that assigns 
to each edge value $0$ (for solid) or value $1$ (for dashed). The graph $A$
is {\it well dashed} if all the $2$-colored $4$-cycles have {\it odd-dashing}, namely
have an odd number of dashed edges. 

\smallskip

An {\it Adinkra} is a well-dashed, $N$-dimensional chromotopology, endowed with a ranking
where all the white vertices have even ranking and all the black vertices have odd ranking.

\smallskip

In the classification of \cite{FaGa05}, the white vertices correspond to the boson fields and their
time derivatives and the black vertices to the fermionic fields and their time derivatives.  
There is an edge between a pair of vertices whenever the corresponding fields are
related by one of the relations \eqref{E:5.2}
 or \eqref{E:5.3}. The edge
is oriented from the white to the black vertex if $\lambda=0$ and viceversa if $\lambda=1$.
The edge is dashed if $c=-1$ and solid if $c=1$.

\smallskip

We refer the reader to { \cite{FaGa05}} and {\cite{DFGHIL07,DFGHIL08,DFGHILM11,DIKLM15,DIKLM16}} for more details about this classification result.

\smallskip

Recall that a binary linear code $L\subset \F_2^N$ is even if every code word $c \in L$ has even weight
$w(c)=\# \{ w_i=1 \}\in 2\mathbf{Z}_{\geq 0}$ and it is doubly even if every code word $c \in L$ has weight that is 
divisible by $4$. Theorem~4.4 of {\cite{DFGHIL08}} shows that Adinkras are classified by doubly-even
binary linear codes, in the sense that every Adinkra $A$ is obtained as a quotient $A=\F_2^N/L$,
for $L$ a doubly-even binary linear code and all the Adinkra structure (that is bipartition, ranking, well-dashing)
is determined by this description. 

\medskip

\begin{proposition} \label{P:5.2}
For a fixed $N\in \N$, consider the set $\cD_N$ of all 
doubly even linear codes $L \subset \F_2^N$,
\[ \mathcal{D}_N:=\{  L\subseteq \F_2^N \,|\, L \text{ doubly even linear code} \}\, . \]
Then, the collection $\{ \mathcal{D}_N \}_{N\geq 1}$ forms a non-unital operad under the composition
\[ \gamma: \mathcal{D}_N \times \mathcal{D}_{k_1}\times \cdots \times \mathcal{D}_{k_N} \to \mathcal{D}_{k_1+\cdots + k_N}, \]
given by
\begin{equation}\label{E:5.4}
\begin{aligned}
 \gamma(L;L_1,\ldots,L_N)=&\left\{ c\in\F_2^{k_1+\cdots + k_N}\,|\,
c=(c^{(1)}c_1,\ldots, c^{(N)} c_N)\, \text{ with } \right. \\
& \left.  (c^{(1)},\ldots,c^{(N)})\in L, \text{ and } c_i \in L_i \right\}\, . 
\end{aligned}
 \end{equation}
\end{proposition}

\medskip

\begin{proof} The code $\gamma(L;L_1,\ldots,L_N)\subset \mathbf{F}_2^{k_1+\cdots + k_N}$
is a linear code since it is a direct sum of copies of the linear codes $L_i$ corresponding
to the subspaces $\F_2^{k_i}$, with $i$ such that $c^{(i)}=1$. 
It is also doubly even since $w(c)=\sum_{i: c^{(i)}=1} w(c_i)$ and each code $L_i$ is
doubly even so each $w(c_i)$ is divisible by $4$, hence $w(c)$ is also divisible by $4$
for all $c\in \gamma(L;L_1,\ldots,L_N)$. 

The composition \eqref{E:5.4}
 satisfies the associativity condition given by the identities
\[ \gamma( \gamma(L^{(N)}; L^{(k_1)}, \ldots, L^{(k_N)}); L^{(r_{1,1})}, \ldots, L^{(r_{1,k_1})}, \ldots, L^{(r_{N,1})}, \ldots, L^{(r_{N,k_N})}) = \]
\[ \gamma(L^{(N)}; \gamma(L^{(k_1)}; L^{(r_{1,1})}, \ldots, L^{(r_{1,k_1})}), \ldots, \gamma(L^{(k_N)}; L^{(r_{N,1})}, \ldots, L^{(r_{N,k_N})})), \]
for $L^{(N)}\in \mathcal{R}_N$, $L^{(k_i)}\in \mathcal{R}_{k_i}$, $i=1,\ldots, N$, 
and $L^{(r_{i,\ell_i})} \in \mathcal{R}_{r_{i,\ell_i}}$ with $\ell_i=1,\ldots, k_i$.
It also satisfies the symmetry conditions given by 
\[ \gamma(\sigma(L); L_{\sigma^{-1}(1)},\ldots,L_{\sigma^{-1}(N)}) =
\tilde\sigma(\gamma(L; L_1,\ldots,L_N)), \]
where on the right-hand-side $\tilde\sigma \in \Sigma_{k_1+\cdots + k_N}$ is the
permutation that splits the set of indices into blocks of $k_i$ indices and permutes the
blocks by $\sigma$, and 
\[ \gamma(L; \sigma_1(L_1),\ldots,\sigma_N(L_N) = \hat\sigma(\gamma (L; L_1,\ldots,L_N)), \]
where on the right-hand-side $\hat\sigma  \in \Sigma_{k_1+\cdots + k_N}$ is the permutation that
acts on the $i$-th block of $k_i$ indices as the permutation $\sigma_i$.

However, the operad is non-unital since $\mathcal{R}_1$ only consists of the trivial subspace of $\mathbf{F}_2$
and that does not satisfy the unital conditions 
$\gamma(L;\mathbf{1},\ldots,\mathbf{1})=L=\gamma(\mathbf{1};L)$. Note that it is
a non-unital operad in the stronger sense of \cite{Markl08}, since the
composition \eqref{E:5.4} is obtained from insertion operations
\[ \circ_i : \mathcal{D}_N \times \mathcal{D}_M \to \mathcal{D}_{N+M-1} \]
of the form
\[ L\circ_i L' =\left\{ c\in \mathbf{F}_2^{N+M-1}\,|\, c=(c^{(1)},\ldots, c^{(i-1)}, c^{(i)} c', c^{(i+1)},\ldots, c^{(N)}) \right. \]
\[ \left. \text{with } (c^{(1)},\ldots, c^{(N)})\in L \text{ and } c'\in L' \right\}\, , \]
with
\[ \gamma(L;L_1,\ldots,L_N)=(\cdots ((L\circ_N L_N)\circ_{N-1} L_{N-1})\cdots \circ_1 L_1)\, . \]
\end{proof}

\medskip

\begin{corollary}  
Let $\mathcal{R}_N$ be the set of representations of the form \eqref{E:5.2} or \eqref{E:5.3} of a
$(1|N)$ SUSY algebra \eqref{E:5.1}. Then, the $\{\mathcal{R}_N \}_{N\geq1}$ form a non-unital operad. 
\end{corollary}

\medskip

\begin{proof} We use the classification of \cite{FaGa05} in terms of Adinkras and the
result of Theorem~4.4 of \cite{DFGHIL08}  recalled above, to 
describe equivalently the set $\mathcal{R}_N$ in terms of doubly even linear codes
in $\mathbf{F}_2^N$. Thus, the operadic composition 
\[ \gamma: \mathcal{R}_N \times \cR_{k_1}\times\cdots\times \mathcal{R}_{k_N} \to \mathcal{R}_{k_1+\cdots+ k_n} \]
associated to a set $(A, A_1,\ldots, A_N)$ of representations, identified with the
corresponding Adinkras $A=\mathbf{F}_2^N/L$, $A_i=\mathbf{F}_2^{k_i}/L_i$, $i=1,\ldots, N$,
with $L$ and the $L_i$ doubly even linear codes, the representation determined
by the Adinkra
\[\gamma(A; A_1,\ldots, A_N):=\F_2^{k_1+\cdots+k_N}/\gamma(L;L_1,\ldots,L_N), \]
where $\gamma(L;L_1,\ldots,L_N)$ is the operadic composition of Proposition~\ref{P:5.2}
\end{proof}

\smallskip

\section{Moufang symmetries of Adinkras}

\bigskip

Hidden symmetries of Adinkras are investigated in this section. It is proved that Adinkras are invariant under Moufang loop symmetries (see~\cite{Sm16} for an introduction to quasigroups,  loops and Moufang loops). These symmetries lead to exploring a very specific class of dessins, corresponding to  the chromotopologies of Adinkras (see~\cite{DIKLM15} for an exposition on the relation between dessins and chromotopologies). 

\smallskip 

In what follows, all our definitions and constructions rely on Sec. 6 of our previous paper \cite{CoMaMar21} concerning Moufang patterns. 

\subsection{Code loops and central extensions}
Sec. II.3, \cite{Che90} presents a generalisation of extension theory defined initially for groups, for the case of loops and quasigroups. This topic has become a popular method, in more recent years~\cite{NaSt08}. We will be considering in what follows the extension of a doubly even code. 

A loop $\mathcal{L}$ is called a code loop of $L$ if $\mathcal {L}$ has a central subgroup $Z$ of order 2 such that $\mathcal{L}/Z\cong L$, as an elementary abelian 2-group and comes equipped with some additional properties. As depicted in Def.2.4~\cite{Hsu} (and in a more general context Prop. 3.1.~\cite{NaSt08}), one needs to define three functions, $\alpha:L\to Z, \phi: L\times L\to Z$ and $\psi: L\times L\times L \to Z$ which play a central role in the definition of the multiplication operation for the code loop and also in the definition of the (central) extension of the (normal) subgroup $Z$ by $L$.  

Consider a doubly even binary code $L\subset \mathbf{F}_2^N$ and let $Z$ be a group of order 2.  Then, by using Thm. 2.5~\cite{Hsu} a code loop of a doubly even code exists and is unique up to isomorphism.  It is important to note that code loops are certain Moufang loop extensions of doubly even binary codes. 

So, supposing that $L$ is a binary linear code in $\F_2^N$ which is doubly even, then one can associate a Moufang loop $\mathcal{L}$ to it such that it satisfies the following exact sequence:

\begin{equation}\label{E:C}1\to Z\to \mathcal{L} \to L\to 1, \end{equation}
 with $L$ a  doubly even binary code, $\mathcal{L}$ the code loop of $L$ (Moufang loop) and $Z$ is a group of order 2 which is a central normal subgroup of $\mathcal{L}$. One has that $L$ is isomorphic to the factor loop $\mathcal{L}/Z$.

\smallskip 

We prove the following statement. 

\medskip 

\begin{theorem}\label{T:6.1.1}
Every Adinkra, interpreted as the quotient $\F_2^N/L$ (with $L$ being a doubly-even binary linear code) 
is invariant under a Moufang loop factor (i.e. a Moufang loop $\mathcal{L}$ quotiented by a finite group $Z$), 
where $Z$ and $\mathcal{L}$ satisfy the following exact sequence:

\[1 \to Z \to \mathcal{L}\to L \to 1,\]

with:

-- $L$ a doubly even binary code

-- $\mathcal{L}$ a Moufang loop 

-- $Z$ a normal subgroup. 
\end{theorem}
\smallskip

\begin{proof}
We apply the discussion above. We have that $L$ is a doubly even binary code, and by Thm. 2.5 \cite{Hsu} we have that a code loop of a doubly even code exists and is unique (up to isomorphism). Using the construction depicted above, we can proceed to defining the following short exact sequence  $$1 \to Z \to \mathcal{L}\to L \to 1,$$ where $L$ is a doubly even binary code, $\mathcal{L}$ is a Moufang loop and $Z$ is a normal subgroup.

Given that Adinkras can be identified to the quotient $\mathbf{F}_2^N/L$ and that $L\cong \mathcal{L}/Z$ this implies now that   $\mathbf{F}_2^N$ is quotiented by a Moufang loop. So, in other words, Adinkras are invariant under Moufang loop symmetries.
\end{proof}
\medskip  
 
 \subsection{Adinkras, graphs and dessins}\label{S:6.2}  
From Sec. \ref{S:5.1} we have that Adinkras form a special class of directed decorated graphs, being $N-$regular, $N-$edged colored bi-partite graphs, where:

\smallskip 

-- vertices are colored black or white and correspond to particles. White vertices correspond to the real bosonic component fields; black vertices correspond to real fermionic component fields.

-- edges correspond to the supersymmetry generators, so that to each vertex $N$ tails are attached.

\medskip 

\begin{proposition}\label{P:6.2.1}
 Let $A$ be the graph of an Adinkra. Then, this graph is invariant under Moufang loop symmetries.
\end{proposition}
\smallskip 

\begin{proof} Previously we proved in Thm.~\ref{T:6.1.1} that every Adinkra is invariant under Moufang loop symmetries. The Moufang loop $\cL$ follows from the exact sequence from formula~\eqref{E:C} and is constructed in a unique way from a given doubly even binary code $L$. So, since the graph reproduces the relations within the Adinkra it follows naturally that the graph $A$ is invariant under Moufang loop symmetries. 
\end{proof}
\medskip

\cite{DIKLM15} implies that, by forgetting all decorations of the Adinkra graphs: dashing, orientation, weights of vertices, (etc) one obtains a coarse version called the chromotopology  of the Adinkra. This object encapsulates the topology of an Adinkra together with the vertex bipartition (coloring each vertex black or white) and the edge $N$-coloring. One can canonically realise these graphs as Grothendieck's dessins d'enfants. The $N$-regular $N$-coloring of edges  gives a cyclic ordering of the edges at each vertex of the graph based on their color, providing thus the structure of a ribbon graph and giving also information concerning the cartographic group of the dessin as we will see in Prop.~\ref{P:6.3.1}. 

\medskip

\subsection{Symmetries of Adinkras' dessins}
In this section the geometries of the dessins, corresponding to the chromotopologies of Adinkras are investigated. The symmetries of those special dessins have repercussions in relation to problems concerning actions of the absolute Galois group $G_{\cQ}$ and relates to ~\cite{CoMaMar20}. This result appears also as a geometric counterpart of the recent result showing that the Grothendieck--Teichm\"uller group $GT$ has dihedral symmetry relations (see~\cite{Co22}).

\smallskip 

From Sec.~\ref{S:6.2} it followed that chromotopologies of Adinkras can be identified to dessins. But dessins can interpreted as being in bijection with transitive subgroups of the symmetric group $\mathbb{S}_n$, generated by 2 elements (for all $n\geq 2$) up to conjugacy (in cases where the dessin corresponds to a degree $n$ map).
This follows from a very geometric construction, which for the sake of clarity, we recall. 

Consider the projective line minus three points $\P^1\setminus \{0,1,\infty\}$. The relation from dessins to degree $n$, finite \'etale coverings of $\P^1\setminus \{0,1,\infty\}$ is well known. Take a point $x\in\P^1\setminus \{0,1,\infty\}$ and its fiber, the set $\{x_1,\cdots,x_n\}.$ Then the topological loops originating from the point $x$ and going around 0, 1 and $\infty$ generate permutations of  the points $\{x_1,\cdots,x_n\}$ denoted respectively $\sigma_0,\sigma_1,\sigma_\infty$
and satisfying $\sigma_0\sigma_1\sigma_\infty=1.$ The permutation $\sigma_1$ is of order 2 and the permutations $\sigma_0,\sigma_1$ generate a subgroup of the symmetric group $\mathbb{S}_n$ which is transitive, in case the covering is connected. 

So, to a dessin we can attach a triplet $ \langle \sigma_0,\sigma_1,\sigma_\infty\rangle \in (\mathbb{S}_n)^3$ where:
 
-- $\sigma_0\sigma_1\sigma_\infty=id_{\mathbb{S}_n}$ 

--  $\langle \sigma_0,\sigma_1,\sigma_\infty\rangle$ is a transitive subgroup of $\mathbb{S}_n.$

In short, denoting by $F_2$ the free group generated by two elements, one can thus identify those objects $\langle \sigma_0,\sigma_1,\sigma_\infty\rangle$ with some group homomorphisms $\tilde{\phi}: F_2 \to \mathbb{S}_n$, for which the group $\tilde{\phi}(F_2)$ acts transitively on $\{1,2,\cdots, n\}$. The equivalence class of these triples forms the dessins; whereas isomorphism classes form the monodromy group of the corresponding dessin. 

\smallskip 

Let us go back to the setting of Sec.\ref{S:5.1} (i.e. $N$-regular bipartite graphs with $2m$ vertices). Let $\dot{\mathcal{D}}$ be a dessin with $n(=Nm)$ edges. Then, we can uniquely associate to it (up to unique isomorphism) a finite set  of size $n$ equipped with a pair of permutations $\dot\sigma, \dot\alpha\in \mathbb{S}_n$  acting transitively on the finite set which consists of the set of edges of $\dot{\mathcal{D}}$ and $\dot\sigma, \dot\alpha$ are the generators of the transitive group corresponding to the Dessin $\dot{\mathcal{D}}$.

\medskip 

\begin{proposition}\label{P:6.3.1}
The dessins corresponding to chromotopologies of Adinkras (being $N$-regular bipartite graphs having a number $n$ of edges and $2m$ of vertices) can be uniquely associate---up to unique isomorphism---to a finite set of size $n$, equipped with a pair of permutations $\dot\sigma, \dot\alpha\in \mathbb{S}_n$, where $\dot\sigma$ and $\dot\alpha$ are both decomposed into $m$ cycles of length $N$.
\end{proposition} 

\smallskip 

\begin{proof} 
Indeed, the graph ought to be an $N$-regular bipartite graph (with $m$ black vertices and $m$ white vertices). The finite set corresponds to the set of (labelled) edges of the chromotopology graph. By definition of the dessin $\dot\sigma$ and $\dot\alpha$ correspond to a decomposition into cycles, where each cycle is a rotation about a black vertex (resp. white) of the finite set. Since we have that the graph is $N$-regular and bipartite with $2m$ vertices we have m blocks of cycles of length $N$ defining $\dot\sigma$ and $\dot\alpha$ respectively. 
\end{proof}
\smallskip 

\begin{proposition} \label{P:6.3.2} 
 Consider the subcategory of dessins  corresponding to Adinkras. Consider a dessin $\dot{\mathcal{D}}$ ($N$-regular bipartite graph on $2m$ vertices). Then, there exists a partition $V_1$ and $V_2$ of the set $\{1,\cdots, Nm\}$ into $m$ subsets of size $N$ defining a specific dessin for which its automorphism group fixing the bipartite blocks (i.e. fixing  the blocks of white and respectively the black vertices set-wise) is $\mathbb{S}_{N2m}$.
 \end{proposition}

\smallskip 

\begin{proof}
 
Consider pair of partitions (say $V_1$ and $V_2$) of the set $\{1,\cdots, mN\}$ into $m$ subsets of size $N$. 
This encodes an $N$-regular bipartite graph. Then, there exists partitions $V_1$ and $V_2$ of $\{1,\cdots, mN\}$ such that:
 the automorphism group fixing the bipartite blocks (i.e. fixing  the blocks of white and respectively the black vertices set-wise) is $\mathbb{S}_{N2m}$. This follows from the Lem 2.3 in~\cite{Ja06}.  
\end{proof}

The following proposition highlights the peculiar symmetries of the Adinkra in relation to the dessins. 
\medskip 

\begin{proposition} \label{P:6.3.3}
 Let us consider an Adinkra $\F_2^N/L$, where $L$ is a doubly even binary code and the corresponding dessin $\dot\cD$ to the Adinkra. Then, $\dot\cD$ is a graph which is invariant under a symmetry group $\dot G$ such that $\dot G \supseteq \cL$, where $\cL$ is the Moufang loop obtained by the extension defined in Equ. \eqref{E:C} of the doubly even binary code $L$.
\end{proposition}
\smallskip 

\begin{proof} Indeed, the dessin corresponding to the Adinkra $\mathbf{F}_2^N/L$  encapsulates the topology of an Adinkra. Somehow, it remains a coarse version of the graph $A$ associated to the Adinkra (because we ``forget" the dashing of edges as well as weights of the vertices). 
In Prop. \ref{P:6.2.1} we showed that the graph $A$ carries the symmetries of the Moufang loop $\mathcal{L}$ associated to its Adinkra. The coarse version (the dessin) therefore inherits the Moufang loop symmetries. However, given that the dashing and weighing of vertices are gone, one can gain more symmetries than in the case of $A$. Therefore, the dessin is invariant under the symmetries of a group $\dot G$ which contains as a subloop the Moufang loop $\mathcal{L}$.
\end{proof}
\smallskip

Note that Adinkras with chromotopology can be associated to an $N$-cube. The $N$-cube has the hypercube symmetry but in particular inherits the dihedral $D_{2N}$ symmetry also. 

\medskip 

\begin{corollary}\label{P:6.3.4}
 Consider an Adinkra $\F_2^N/L$. Then, it has the structure of the quotient of an $N$-hypercube by a Moufang loop and inherits dihedral symmetries, quotiented by relations given by $L$. 
\end{corollary}

\smallskip 

\begin{proof} 
The proof follows from the fact that every connected Adinkra chromotopology is isomorphic to a quotient of a colored $N$-dimensional cube by the code of the chromotopology and from the natural symmetries of the hypercube. 

Finally,  since the $N$-hypercube has the symmetry of the wreath product of $\mathbb{S}_2\, wr\, \mathbb{S}_N$ and $\mathbb{S}_N$ has as a subgroup the dihedral group $D_{2N}$, we can say that this $N$-cube  inherits the symmetries of the dihedral $D_{2N}$ group. So, the Adinkra (and chromotopology) have the symmetry of a dihedral $D_{2N}$ group quotiented by the relations given by $L$.  
\end{proof}

\bigskip 

\thanks{{\bf Acknowledgements.} N. C. Combe acknowledges support from the Minerva fast track grant from the Max Planck Institute for Mathematics in the Sciences, in Leipzig.
Yu. Manin acknowledges the continuing strong support from the Max Planck Institute for Mathematics in Bonn.
M. Marcolli acknowledges support from NSF grant DMS- 2104330.} 

\bigskip

\bibliographystyle{amsalpha}

\begin{thebibliography}{A}


\bibitem[BoMa]{BoMa08}
  D.~Borisov, Yu. Manin, {\it Generalized operads and their inner cohomomorhisms,}
 In: Geometry and Dynamics of Groups
and Spaces (In memory of Aleksander Reznikov). Ed. by M. Kapranov et al.
Progress in Math., vol. 265 (2008), 
Birkh\"auser, Boston, pp. 247--308.
arXiv:math.CT/0609748.


\smallskip

\bibitem[BrHePo] {BrHePo20}
 U.  Bruzzo,  D. Hern\'andez Ruip\'erez,  A.  Polischuk, {\it Notes on fundamental
algebraic supergeometry. Hilbert and Picard superschemes,} arXiv:2008.00700 [math.AG]. 

\smallskip

\bibitem[Co] {Co22} 
N. Combe, {\it Dihedral symmetries of the Grothendieck--Teichm\"uller group,} To appear. 

\smallskip

\bibitem[Che]{Che90}
 O. Chein, {\it Examples and methods of construction,} Chapter II in Quasigroups and Loops:
Theory and Applications (O. Chein, H. O. Pflugfelder, J. D. H. Smith, eds.). Heldermann
Verlag, Berlin (1990), pp. 27--93.

\smallskip

\bibitem[ChCorGi] {ChCorGi19} 
C.  Chenavier,  C.  Cordero,  S.  Giraudo,  {\it Quotients of the magmatic
operad: lattice structures and convergent rewrite systems,} arXiv:1809.05083v2.

\smallskip

\bibitem[CoMaMar1]{CoMaMar20}
 N. Combe, Yu. Manin., M. Marcolli, {\it Dessins for modular operad and Grothendieck--Teichm\"uler group,} ``Topology and Geometry A Collection of Essays Dedicated to Vladimir G. Turaev'', European Mathematical Society (2021), arXiv:math.AG/2006.13663.


\smallskip



\bibitem[CoMaMar2]{CoMaMar21}
 N. C. Combe, Yu. Manin, M. Marcolli, {\it Moufang patterns and geometry of information,}
To be published in the Collection, dedicated to Don Zagier, in Pure and Applied Math. Quarterly.
arXiv:2107.07486, 42 pp.

\smallskip

\bibitem[CoMaMar3]{CoMaMar22}
 N.  C. Combe, Yu. Manin, M. Marcolli, {\it Quantum operads,}
to be published in the ``C. N.  Yang at 100'' volume, arXiv:2112.15237, 34 pp.

\smallskip 

\bibitem[De]{De87}
 P. Deligne. {\it Lettre \`a Manin,} Princeton. 
URL: https:// publications.ias.edu/

sites/default/files/lettre-a-manin-1987-09-25.pdf


\smallskip

\bibitem[DFGHIL1]{DFGHIL07} C.F.~Doran, M.G.~Faux, S.J.~Gates, Jr., T.~H\"ubsch, K.M.~Iga, G.D.~Landweber, 
{\it On graph-theoretic identifications of Adinkras, supersymmetry representations and superfields}, 
Int. J. Mod. Phys. A 22 (2007), pp. 869--930.

\smallskip

\bibitem[DFGHIL2]{DFGHIL08} C.F.~Doran,  M.G.~Faux,  S.J.~Gates, Jr., T.~H\"ubsch, K.M.~Iga, G.D.~Landweber, 
{\it Relating doubly-even error-correcting codes, graphs, and irreducible representations of $N$-supersymmetry,} in ``Discrete and computational mathematics", pp.53--71, Nova Sci. Publ., (2008).

\smallskip

\bibitem[DFGHILM]{DFGHILM11} C.F.~Doran,  M.G.~Faux,  S.J.~Gates, Jr., T.~H\"ubsch, K.M.~Iga, G.D.~Landweber, R.L.~Miller, {\it Codes and supersymmetry in one dimension}, Adv. Theor. Math. Phys., Vol.15 (2011) 
N.6, pp. 1909--1970.

\smallskip

\bibitem[DIKLM1]{DIKLM15} C.F.~Doran, K.~Iga, J.~Kostiuk, G.~Landweber, S.~M\'endez-Diez, {\it Geometrization
of $N$-extended $1$-dimensional supersymmetry algebras, I}, 
Adv. Theor. Math. Phys. 19 (2015), no. 5,  pp. 1043--1113. 

\smallskip

\bibitem[DIKLM2]{DIKLM16} C.F.~Doran, K.~Iga, J.~Kostiuk, S.~M\'endez-Diez,
 {\it Geometrization of $N$-extended $1$-dimensional supersymmetry algebras, II}, arXiv:1610.09983.

\smallskip 

\bibitem[FaGa]{FaGa05} M.~Faux, S.J.~Gates, Jr., {\it Adinkras: a graphical technology for supersymmetric representation theory}, Phys. Rev. D 71(3) (2005), 065002.

\smallskip 

\bibitem[FeKaPo]{FeKaPo20} 
G. Felder,  D. Kazhdan, A.  Polishchuk,  {\it The moduli space 
of stable supercurves and its canonical line bundle,} arXiv:2006.13271v2 [math.AG].

\smallskip

\bibitem[GeKa]{GeKa98}
 E. Getzler, M.  Kapranov, {\it Modular operads,} In: Compositio Math. (1998),
110:1, pp. 65--125.

\smallskip 

\bibitem[HoKrStT]{HoKrStT11}
 H.  Hohnhold,  M.  Kreck,  St. Stolz, P. Teichner, 
{\it Differential forms and $0$--dimensional supersymmetric field theories,}
Quantum Topology  2 (2011), pp. 1--41.

\smallskip


\bibitem[Hsu]{Hsu}
T. Hsu, {\it Explicit constructions of code loops as centrally twisted products,}
 Math. Proc. Camb. Phil. Soc. 128 (2000), pp. 223--232

\smallskip

\bibitem[Ja]{Ja06}
 J.P. James, {\it Partition actions of symmetric groups and regular bipartite graphs,}
 Bull. London Math. Soc. 38 (2006), pp. 224–232

\smallskip

\bibitem[KaSch]{KaSh06}
 M.  Kashiwara,  P.  Schapira,  {\it Categories and sheaves,}  Springer Verlag (2006),
pp. x + 497.

\smallskip

\bibitem[KeMaWu]{KeMaWu22}
E. Kessler, Yu. Manin,  Y. Wu, {\it Moduli spaces of SUSY curves
and their operads,} arXiv:2202.10321 [mathAG], 21 pp.

\smallskip

\bibitem[Ke]{Ke01}
 B. Keller, {\it Introduction to $A_{\infty}$--algebras and modules,} In "Homology,
Homotopy and Applications" 3 (2001), pp.1--35 (electronic), arXiv:math/9910179v2, 31 pp.

\bibitem[Ma]{Ma88}  Yu.  Manin,  {\it Quantum groups and non--commutative geometry,} Publ. de CRM,
Universit\'e  de Montr\'eal (1988), 91 pp.

\smallskip

\bibitem[MaPe]{MaPe89}
 Yu. Manin, I. Penkov, {\it The formalism of left and right connections
on supermanifolds,} In "Lectures on Supermanifolds, Geometrical Methods
and Conformal Groups", World Scientific (1989), pp. 3--13.


\smallskip

\bibitem[MaVa]{MaVa20}
 Yu.  Manin,  B.  Vallette, {\it Monoidal structures on the categories of 
quadratic data,} Documenta Mathematica 25 (2020), pp. 1653 -- 1712,
arXiv:1902.03778.

\smallskip

\bibitem[Markl]{Markl08} M.~Markl, {\it Operads and PROPs}, 
Handbook of algebra. Vol. 5, 87–140, Elsevier/North-Holland, 2008.

\smallskip

\bibitem[NaSt]{NaSt08} P. T. Nagy, K. Strambach, 
{\it Schreier loops,} Czechoslovak Mathematical Journal, Vol. 58 (2008), No. 3, 759--786

\smallskip

\bibitem[Sm]{Sm16}
 J. D. H. Smith, {\it Quantum quasigroups and loops,}  Journ. of Algebra 456 (2016),
pp. 46--75.

\smallskip

\bibitem[Val]{Val04}
 B.  Vallette, {\it A Koszul duality for PROPs,} arXiv:math0411542v3, 78 pp.

\smallskip

\bibitem[WuYau]{WuYau06}
Y.  Wu, Sh.-T. Yau, {\it Comparison theorems of phylogenetic spaces
and the moduli spaces of curves,} arXiv:2006.04319 [math.AG].

\end{thebibliography}

\medskip

{\bf Noemie C. Combe, Max Planck Institute for Mathematics in the Sciences, Inselstrasse 22, 04103 Leipzig, Germany. \\

Yuri I. Manin, Max Planck Institute for Mathematics, Vivatsgasse 7, 53111 Bonn. \\

Matilde Marcolli, Math. Department, Mail Code 253-37, Caltech, 1200 E.California Blvd., Pasadena, CA 91125, USA.}

\end{document}